\documentstyle[11pt]{article}
 \newcommand{\qb}[2]{{\left [{#1 \atop #2} \right]}}
 \newlength{\standardunitlength}
\newtheorem{cor}{Corollary} \newtheorem{lemma}{Lemma}
\newtheorem{theorem}{Theorem} \newtheorem{prop}{Proposition}
\newenvironment{proof}{\noindent {\sc Proof:}}{$\Box$ \vspace{2 ex}}

\begin{document}

\begin{center} Applications of Symmetric Functions to Cycle and
Increasing Subsequence Structure after Shuffles (Part 2) \end{center}

\begin{center}
By Jason Fulman
\end{center}

\begin{center}
Version 1: March 31, 2001
\end{center}

\begin{center}
Current version : April 14, 2001
\end{center}

\begin{center}
Stanford University
\end{center}

\begin{center}
Department of Mathematics
\end{center}

\begin{center}
Building 380, MC 2125
\end{center}

\begin{center}
Stanford, CA 94305
\end{center}

\begin{center}
email:fulman@math.stanford.edu
\end{center}

\newpage \begin{abstract} Using the Berele/Remmel/Kerov/Vershik
variation of the Robinson-Schensted-Knuth correspondence, we study the
cycle and increasing subsequence structure after various methods of
shuffling. One consequence is a cycle index for shuffles like: cut a
deck into two roughly equal piles, thoroughly mix the first pile and
then riffle it with the second pile. Conclusions are drawn concerning
the distribution of fixed points and the asymptotic distribution of
cycle structure. An upper bound on the convergence rate is
given. Connections are made with extended Schur functions and with
point process work of Baik and Rains.  \end{abstract}


\newpage

\section{Introduction} \label{Introduction}

	In an unpublished effort to study the way real people shuffle
cards, Gilbert-Shannon-Reeds introduced the following model, called
$k$-riffle shuffling. Given a deck of $n$ cards, one cuts it into $k$
piles with probability of pile sizes $j_1,\cdots,j_k$ given by
$\frac{{n \choose j_1,\cdots,j_k}}{k^n}$. Then cards are dropped from
the packets with probability proportional to the pile size at a given
time (thus if the current pile sizes are $A_1,\cdots,A_k$, the next
card is dropped from pile $i$ with probability
$\frac{A_i}{A_1+\cdots+A_k}$). A celebrated result of Bayer and
Diaconis \cite{BD} is that $\frac{3}{2}log_2(n)$ two-shuffles are
necessary and suffice for randomness.

	One of the most remarkable properties of GSR $k$-shuffles is
the following. Since $k$-shuffles induce a probability measure on
conjugacy classes of $S_n$, they induce a probability measure on
partitions of $n$. Consider the factorization of a random degree $n$
polynomial over a finite field $F_q$ into irreducibles. The degrees of
the irreducible factors of a randomly chosen degree $n$ polynomial
also give a random partition of $n$. The fundamental result of
Diaconis-McGrath-Pitman (DMP) \cite{DMP} is that this measure on
partitions of $n$ agrees with the measure induced by card shuffling
when $k=q$. This allowed natural questions on shuffling to be reduced
to known results on factors of polynomials and vice versa.

	The DMP theorem is deep and has connections to many parts of
mathematics (e.g. Hochschild homology \cite{H}, dynamical systems
\cite{La}, and Lie theory \cite{F1},\cite{F2},\cite{F3}). Stanley
\cite{Sta} gave a proof of the DMP theorem using ideas from symmetric
function theory. He also related the Robinson-Schensted-Knuth shape of
a permutation after a shuffle to Schur functions and gave connections
with work of the random matrix community. Following the appearance of
\cite{Sta}, the paper \cite{FDeal} gave a different symmetric function
theoretic proof of the DMP result. Although more complicated than
Stanley's proof (in particular it needed the RSK algorithm), it had
the merit of suggesting natural extensions of the DMP result.

	This note is a continuation of \cite{FDeal} (Part 1) and was
motivated by an effort to understand the relation between it and the
paper \cite{KV}, in particular the fact that the $S_{\lambda}$ used in
Part 1 are extended Schur functions. Section \ref{Extended} defines
models of card shuffling called $(\vec{\alpha},\vec{\beta},\gamma)$
shuffles which include the GSR shuffles. This model contains other
shuffles of interest such as: given a deck of $n$ cards, cut off
binomial(n,1/2) many cards as in the GSR $2$-shuffle, shuffle them
thoroughly, and then riffle them with the remaining cards (this
special case was first studied in \cite{DFP}). It is proved that if
one applies the RSK correspondence to a permutation distributed as a
$(\vec{\alpha},\vec{\beta},\gamma)$ shuffle, then the probability of
getting any recording tableau of shape $\lambda$ is the extended Schur
function $S_{\lambda}(\vec{\alpha},\vec{\beta},\gamma)$. For the case
$\gamma=0$ this result will be reduced to Proposition 3 of the paper
\cite{KV}. Since Proposition 3 of \cite{KV} is incorrect for $\gamma
\neq 0$, Section \ref{Extended} undertakes the repair needed for the
applications to card shuffling. Section \ref{Cauchy} proves
Cauchy-type identities for the extended Schur functions. (The method
of proof closely follows that of \cite{TW} for Schur functions; since
we need the Cauchy type identities and similar reasoning later in the
paper, we include the details). Section \ref{Connect} connects the
shuffling models with work of Baik and Rains \cite{BaRa} on increasing
subsequences for point process models. The final section applies the
work of earlier sections to find formulas for cycle structure after
$(\vec{\alpha},\vec{\beta},\gamma)$ shuffles; for instance it is
proved that after such a shuffle on a deck of size $n$, the expected
number of fixed points is the sum of the first $n$ extended power sum
symmetric functions evaluated at the relevant parameters. An upper
bound on the convergence rate of the shuffles is derived.

	As Stanley \cite{Sta} notes, the results about longest
increasing subsequences after riffle shuffles (which is the same as
longest weakly increasing subsequences in random words) tie in with
work of the random matrix community
(e.g. \cite{TW},\cite{K},\cite{Jo}). We expect that the connection
between card shuffling and random matrix theory holds for all of the
shuffles considered here (and possibly for the shuffles in
\cite{F2},\cite{F3}). Some work has been done in this direction (not
in the language of card shuffling) and is nicely surveyed in
\cite{BorO}. It is remarkable that the unified viewpoint of
\cite{BorO} uses the represenation theory of the infinite symmetric
group (extended Schur functions), which are also the parameterizing
data in this paper.

\section{Extended Schur Functions and the RSK Correspondence}
\label{Extended}

	The extended complete symmetric functions
$\tilde{h}_k(\alpha,\beta,\gamma)$ are defined by the generating
function \[ \sum_{k=0}^{\infty} \tilde{h}_k(\alpha,\beta,\gamma) z^k =
e^{\gamma z} \prod_{i \geq 1} \frac{1+\beta_i z}{1-\alpha_i z}.\] For
$\lambda=(\lambda_1,\cdots,\lambda_n)$, the extended Schur functions
are defined by \[ \tilde{s}_{\lambda} =
det(\tilde{h}_{\lambda_i-i+j})_{i,j=1}^n .\] Since the extended Schur
functions give the characters of the infinite symmetric group, they
are very natural objects.

	Thoma proves the following result.

\begin{theorem} \label{Thoma's} (\cite{T}) Let $G(z) =
\sum_{k=0}^{\infty} g_k z^k$ be such that $g_0=1$ and all $g_k \geq
0$. Then

\[ det(g_{\lambda_i-i+j})_{i,j=1}^n \geq 0 \] for all partitions
$\lambda$ if and only if

\[ G(z) = e^{\gamma z} \prod_{i \geq 1} \frac{1+\beta_i z}{1-\alpha_i
z} \] where $\gamma \geq 0$ and $\sum \beta_i, \sum \alpha_i$ are
convergent series of positive numbers. \end{theorem}

	Given Thoma's result, it is natural to interpret the extended
Schur functions as probabilities. This is the topic of the paper
\cite{KV}; however their Proposition 3 is false for $\gamma \neq
0$. This section repairs their Proposition 3 and makes a connection
with card shuffling.

	We suppose that $\gamma+\sum \alpha_i+\sum \beta_i=1$ and that
$\gamma \geq 0$, $\alpha_i,\beta_i \geq 0$ for all $i$. Using these
parameters, we define a random permutation on $n$ symbols as
follows. First, create a word of length $n$ by choosing letters $n$
times independently according to the rule that one picks $i>0$ with
probability $\alpha_i$, $i<0$ with probability $\beta_i$, and $i=0$
with probability $\gamma$. We use the usual ordering $\cdots < -1 < 0
< 1 < \cdots$ on the integers. Starting with the smallest negative
symbol which appears in the word, let $m$ be the number of times it
appears. Then write $\{1,2,\cdots,m\}$ under its appearances in {\it
decreasing} order from left to write. If the next negative symbol
appears $k$ times write $\{m+1,\cdots,m+k\}$ under its appearances,
again in decreasing order from left to write. After finishing with the
negative symbols, proceed to the $0$'s. Letting $r$ be the number of
$0$'s, choose a random permutation of the relevant $r$ consecutive
integers and write it under the $0$'s. Finally, move to the positive
symbols. Supposing that the smallest positive symbol appears $s$
times, write the relevant $s$ consecutive integers under its
appearances in {\it increasing} order from left to right.

 	The best way to understand this procedure is through an
example. Given the string \[ -2 \ 0 \ 1 \ 0 \ 0 \ 2 \ -1 \ -2 \ -1 \ 1
\] one obtains each of the six permutations

\[ 2 \ 5 \ 8 \ 6 \ 7 \ 10 \ 4 \ 1 \ 3 \ 9 \] \[ 2 \ 5 \ 8 \ 7 \ 6 \ 10
\ 4 \ 1 \ 3 \ 9 \] \[ 2 \ 6 \ 8 \ 5 \ 7 \ 10 \ 4 \ 1 \ 3 \ 9 \] \[ 2 \
6 \ 8 \ 7 \ 5 \ 10 \ 4 \ 1 \ 3 \ 9 \] \[ 2 \ 7 \ 8 \ 5 \ 6 \ 10 \ 4 \
1 \ 3 \ 9 \] \[ 2 \ 7 \ 8 \ 6 \ 5 \ 10 \ 4 \ 1 \ 3 \ 9 \] with
probability $1/6$. In all cases the $1,2$ correspond to the $-2$'s,
the $3,4$ correspond to the $-1$'s, the $8,9$ correspond to the $1$'s
and the $10$ corresponds to the $2$. The symbols $5,6,7$ correspond to
the $0$'s and there are six possible permutations of these symbols. We
call this probability measure on permutations a
$(\vec{\alpha},\vec{\beta},\gamma)$ shuffle.

	The following elementary result (generalizing results in
\cite{BD} and \cite{DFP}) gives physical descriptions of these
shuffles and explains how they convolve. The proof method follows that
of \cite{BD}.

\begin{prop} \label{describe}
\begin{enumerate}

\item A $(\vec{\alpha},\vec{\beta},\gamma)$ shuffle is equivalent to
the following procedure. Cut the $n$ card deck into piles with sizes
$X_i$ indexed by the integers, where the probability of having
$X_i=x_i$ for all $i$ is equal to \[ \frac{n!}{\prod_{i=-
\infty}^{\infty} x_i!}  \gamma^{x_0} \prod_{i >0} \alpha_i^{x_i}
\prod_{i<0} \beta_i^{x_i}. \] The top cards go to the non-empty pile
with smallest index, the next batch of cards goes to the pile with
second smallest index, and so on. Then mix the pile indexed by $0$
until it is a random permutation, and turn upside down all of the
piles with negative indices. Finally, riffle the piles together as in
the first paragraph of the introduction and look at the underlying
permutation (i.e. ignore the fact that some cards are upside down).

\item The inverse of a $(\vec{\alpha},\vec{\beta},\gamma)$ shuffle is
equivalent to the following procedure. Randomly label each card of the
deck, picking label $0$ with probability $\gamma$, label $i>0$ with
probability $\alpha_i$ and label $i<0$ with probability
$\beta_i$. Deal cards into piles indexed by the labels, where cards
with negative or zero label are dealt face down and cards with
positive label are dealt face up. Then mix the pile labelled $0$ so
that it is a random permutation and turn all of the face up piles face
down. Finally pick up the piles by keeping piles with smaller labels
on top.

\item Performing a $(\vec{\alpha},\vec{\beta},\gamma)$ shuffle $k$
times is the same as performing the following shuffle. One cuts into
piles with labels given by $k$-tuples of integers $(z_1,\cdots,z_k)$
ordered according to the following rule:

\begin{enumerate}

\item $(z_1,\cdots,z_k) < (z_1',\cdots,z_k')$ if $z_1<z_1'$.

\item $(z_1,\cdots,z_k) < (z_1',\cdots,z_k')$ if $z_1=z_1' \geq 0$ and
$(z_2,\cdots,z_k) < (z_2',\cdots,z_k')$.

\item $(z_1,\cdots,z_k) < (z_1',\cdots,z_k')$ if $z_1=z_1' < 0$ and
 $(z_2,\cdots,z_k) > (z_2',\cdots,z_k')$.

\end{enumerate} The pile is assigned probability equal to the product
of the probabilities of the symbols in the $k$ tuple. Then the shuffle
proceeds as in part 1, where negative piles (piles where the product
of the coordinates of the $k$ tuple are negative) are turned upside
down and piles with some coordinate equal to 0 are perfectly mixed
before the piles are all riffled together.  \end{enumerate} \end{prop}

{\bf Examples} As an example of Proposition \ref{describe}, consider
an $(\alpha_1,\alpha_2; \beta_1,\beta_2;\gamma)$ shuffle with
$n=11$. For part 1, it may turn out that $X_{-2}=2$, $X_{-1}=1$,
$X_0=3$, $X_1=2$, and $X_2=3$. Then the deck is cut into piles
$\{1,2\}$, $\{3\}$, $\{4,5,6\}$, $\{7,8\}$, $\{9,10,11\}$. The first
two piles are turned upside down and the third pile is completely
randomized, which might yield piles $\{2,1\}$, $\{3\}$, $\{5,4,6\}$,
$\{7,8\}$, $\{9,10,11\}$. Then these piles are riffled together as in
the GSR shuffle. This might yield the permutation

\[ 5 \ 2 \ 7 \ 4 \ 8 \ 9 \ 10 \ 3 \ 1 \ 11 \ 6 .\] The inverse
description (part 2) would amount to labelling cards 2,9 with $-2$,
card 8 with $-1$, cards 1,4,11 with $0$, card 3,5 with $1$, and cards
6,7,10 with $2$, and then mixing the $0$ pile as $4,1,11$. Note that
this leads to the permutation (inverse to the previous permutation)

\[ 9 \ 2 \ 8 \ 4 \ 1 \ 11 \ 3 \ 5 \ 6 \ 7 \ 10 .\]

	As an example of part 3, note that doing a
$(\alpha_1;\beta_1;0)$ shuffle twice does not give a
$(\vec{\alpha},\vec{\beta},\gamma)$ shuffle, but rather gives a
shuffle with 4 piles in the order $(-1,1),(-1,-1),(1,-1),(1,1)$ where
pile 1 has probability $\beta_1 \alpha_1$, pile 2 has probability
$\beta_1 \beta_1$, pile 3 has probability $\alpha_1 \beta_1$ and pile
4 has probability $\alpha_1 \alpha_1$. Piles 1 and 3 are turned upside
down before the riffling takes place. From part 1 of this paper one
can still analyze the cycle structure and RSK shape of these shuffles
even though they aren't $(\vec{\alpha},\vec{\beta},\gamma)$
shuffles. (Actually part 1 of this paper looks at shuffles conjugate
to these shuffles by the longest element; this clearly has no effect
on the cycle index and has no effect on the RSK shape by a result of
Sch$\ddot{u}$tzenberger exposited as Theorem A1.2.10 in \cite{Sta2}).

	As another example of part 3, note that a shuffle with
parameters $(\alpha_1;0;\gamma)$ repeated twice gives a shuffle with 4
piles in the order $(0,0),(0,1),(1,0),(1,1)$ where the first 3 piles
are completely mixed before all piles are riffled together. This is
clearly the same as a $(\alpha_1^2;0;1-\alpha_1^2)$ shuffle, agreeing
with Lemma 2.1 of \cite{DFP}.

	Berele and Remmel \cite{BR} and independently Kerov and
Vershik \cite{KV} consider the following analog of the RSK
correspondence. Given a word on the symbols $\{\pm 1,\pm 2,\cdots\}$
one runs the RSK correspondence with the ammendments that a negative
symbol is allowed to bump itself, but that a positive symbol can't
bump itself. For example the word \[ 1 \ -1 \ 2 \ -2 \ 1 \ 1 \ -2 \]
has insertion tableau $P$ and recording tableau $Q$ respectively equal
to

\[ \begin{array}{c c c c c c c c c c c}
             & & & -2 & 1 & 1 &  &   &  &  & \\
            &&  &  -2 & 2 &  &  &&& &\\
             &&&    -1 &  &  &  & &&&\\
             &&& 1 &&&&&&& \end{array}  \]

\[ \begin{array}{c c c c c c c c c c c}
             & & & 1 & 3 & 6 &  &   &  &   & \\
            &&   &  2 & 5 &  & && & &\\
             &&&    4 &  &  &  & &&&\\
              &&& 7 &&&&&&& \end{array}  \]

\begin{theorem} \label{Ber1} (\cite{BR},\cite{KV}) The above variation
on the Robinson-Schensted-Knuth correspondence gives a bijection
between words of length $n$ from the alphabet of integers with the
symbol $i$ appearing $n_i$ times and pairs $(P,Q)$ where

\begin{enumerate}
\item The symbol $i$ occurs $n_i$ times in $P$.
\item The entries of $P$ are weakly increasing in rows and columns.
\item Each positive symbol occurs at most once in each column of $P$
and each negative symbol occurs at most once in each row of $P$.
\item $Q$ is a standard Young tableau on the symbols $\{1,\cdots,n\}$.
\end{enumerate} Furthermore, \[ S_{\lambda}(\vec{\alpha},\vec{\beta},0) = \sum_{P \atop shape(P)=\lambda} \prod_{i > 0} \alpha_i^{n_i(P)} \prod_{i<0} \beta_i^{n_i(P)}.\] \end{theorem}

	Theorem \ref{probinter} and Corollary \ref{c1} connect card
shuffling to the extended Schur functions. Related results can be
found \cite{Sta} and \cite{FDeal} (in particular, Corollary \ref{c1}
was proved in \cite{Sta} for usual Schur functions).

\begin{theorem} \label{probinter} Let $\pi$ be distributed as a
permutation under a $(\vec{\alpha},\vec{\beta},\gamma)$ shuffle. Let $Q$
be any standard Young tableaux of shape $\lambda$. Then the
probability that $\pi$ has Robinson-Schensted-Knuth recording tableau
equal to $Q$ is
$S_{\lambda}(\vec{\alpha},\vec{\beta},\gamma)$. \end{theorem}

\begin{proof} First suppose that $\gamma = 0$. As indicated earlier in
this section, each length $n$ word $w$ on the symbols $\{\pm 1, \pm 2,
\cdots \}$ defines a permutation $\pi$. From this construction, it is
easy to see that the recording tableau of $w$ under the BRKV variation
of the RSK algorithm is equal to the recording tableau of $\pi$ under
the RSK algorithm. Thus it is enough to prove that the probability
that the word $w$ has BRKV recording tableau $Q$ is
$S_{\lambda}(\vec{\alpha},\vec{\beta},0)$. This is immediate from
Theorem \ref{Ber1}.
	
	Now the case $\gamma \neq 0$ can be handled by introducing $m$
extra symbols between $0$ and $1$--call them
$1/(m+1),2/(m+1),\cdots,m/(m+1)$ and choosing each with probability
$\gamma / m$. Thus the random word is on $\{\pm 1,\pm 2, \cdots\}$ and
these extra symbols. Each word defines exactly one permutation--the
symbols $1/(m+1),2/(m+1),\cdots,m/(m+1)$ are treated as positive. By
the previous paragraph, the probabiilty of obtaining recording tableau
$Q$ is equal to $S_{\lambda}(\vec{\alpha},\vec{\beta})$ where the
associated $\tilde{h}_{k}$ are defined by \[ \sum_{k=0}^{\infty}
\tilde{h}_k(\alpha,\beta) z^k = (\frac{1}{1-\gamma z /m})^m \prod_{i
\geq 1} \frac{1+\beta_i z}{1-\alpha_i z}.\] As $m \rightarrow \infty$,
this distribution on permutations converges to that of a
$(\vec{\alpha},\vec{\beta},\gamma)$ shuffle, and the generating
function of the $\tilde{h}_k$ converges to \[ \sum_{k=0}^{\infty}
\tilde{h}_k(\alpha,\beta,\gamma) z^k = e^{\gamma z} \prod_{i \geq 1}
\frac{1+\beta_i z}{1-\alpha_i z}.\] \end{proof}

\begin{cor} \label{c1} Let $f_{\lambda}$ be the number of standard
Young tableau of shape $\lambda$. Let $\pi$ be distributed as a
permutation under a $(\vec{\alpha},\vec{\beta},\gamma)$ shuffle. Then
the probability that $\pi$ has Robinson-Schensted-Knuth shape $\lambda$
is equal to $f_{\lambda}
S_{\lambda}(\vec{\alpha},\vec{\beta},\gamma)$. \end{cor}

\section{Gessel's Theorem, Szego's Theorem, and Cauchy Identities}
\label{Cauchy}

	Letting $f(z)= \sum_{-\infty}^{\infty} c_k z^k$ be a function
on the unit circle, the Toeplitz determinant $D_{n}(f)$ is defined as
the determinant of the $n \times n$ Toeplitz matrix

\[ \left( \begin{array}{c c c c}
		c_0 & c_1 & \cdots &c_{n-1} \\
		c_{-1} & c_0 & \cdots & c_{n-2}\\
		\cdots & \cdots & \cdots & \cdots \\
		c_{-(n-1)} & c_{-(n-2)} & \cdots & c_0
	  \end{array} \right). \] The function $f$ is called the
symbol of the Toeplitz determinant. The strong Szego theorem for
Toeplitz determinants \cite{BS} states that under mild conditions on
the symbol, \[ D_{n}(e^{ \sum_{-\infty}^{\infty} c_k z^k}) \sim
e^{nc_0 + \sum_{k=1}^{\infty} kc_kc_{-k}}.\]

	The aim of this section is to prove the following two
identities:

 \[ \sum_{\lambda \atop l(\lambda)
\leq n} S_{\lambda}(\alpha,\beta,\gamma) s_{\lambda}(x) =
D_{n} \left(e^{\gamma z} \prod_{r=1}^{\infty} \frac{1+\beta_r
z}{(1-x_r/z)(1-\alpha_r z)}\right) \]

\[ \sum_{\lambda} S_{\lambda}(\alpha,\beta,\gamma) s_{\lambda}(x) =
\sum_{\lambda} \frac{1}{z_{\lambda}} \tilde{p}_{\lambda}
(\alpha,\beta,\gamma) p_{\lambda}(x) \]

	The first identity will be proved using a technique of Gessel
\cite{G}, and the second identity will follow from the first by
applying the strong Szego limit theorem for Toeplitz determinants (an
idea used in Tracy and Widom \cite{TW}).

\begin{theorem} \label{likeGessel}

\[ \sum_{\lambda \atop l(\lambda) \leq n}
S_{\lambda}(\alpha,\beta,\gamma) s_{\lambda}(x) = D_{n}
\left(e^{\gamma z} \prod_{r=1}^{\infty} \frac{1+\beta_r
z}{(1-x_r/z)(1-\alpha_r z)} \right) \] \end{theorem}

\begin{proof} Let $A$ be the $n \times \infty$ matrix
$(\tilde{h}_{j-i}(x))$ ($i \geq 1, 1 \leq j \leq n$). Let $B$ be the
$\infty \times n$ matrix $(h_{i-j}(x))$ ($i \geq 1, 1 \leq j \leq
n$). We evaluate in two ways $det(AB)$. 

	One on hand, 

\[ det(AB) = det \left(\sum_{k=0}^{\infty}
\tilde{h}_{k-i}(\alpha,\beta,\gamma) h_{k-j}(x)\right).\] Clearly this is a
Toeplitz determinant. The symbol is

\begin{eqnarray*}
\sum_{j=-\infty}^{\infty} z^j \sum_{k=0}^{\infty} \tilde{h}_k(\alpha,\beta,\gamma) h_{k-j}(x) & = & \sum_{k=0}^{\infty} \tilde{h}_k(\alpha,\beta,\gamma) \sum_{j=-\infty}^{\infty} z^j h_{k-j}(x)\\
& = &  \sum_{k=0}^{\infty} z^k \tilde{h}_k(\alpha,\beta,\gamma) \sum_{j=-\infty}^{\infty} z^{j-k} h_{k-j}(x)\\
& = & e^{\gamma z} \prod_{r=1}^{\infty} \frac{1+\beta_r
z}{(1-x_r/z)(1-\alpha_r z)}.
\end{eqnarray*}

	On the other hand, the Cauchy-Binet formula gives that

\[ det(AB) = \sum_{S=s_1<\cdots<s_n} det(A|S) det(S|B) \] where
$(A|S)$ (resp. $(B|S)$) is the $n \times n$ matrix formed by using the
columns (resp. rows) indexed by $S$. Writing $S=\{\lambda_{n+1-i}+i\}$
associates the subsets $S$ with partitions with at most $n$
parts. Then
\[ \sum_{S=s_1<\cdots<s_n} det(A|S) det(S|B) = \sum_{\lambda \atop
l(\lambda) \leq n} S_{\lambda}(\alpha,\beta,\gamma) s_{\lambda}(x) \]
as desired.
\end{proof}

	We also use the well known Polya cycle index.

\begin{lemma} \label{Polya} Let $m_i(\lambda)$ be the number of parts
of size $i$ in $\lambda$. Then \[ \sum_{\lambda} \frac{1}{z_{\lambda}}
\prod_i x_i^{m_i(\lambda)} = e^{\sum_{k \geq 1} \frac{x_k}{k}}.\]
\end{lemma}

\begin{proof} The coefficient of $\prod_i x_i^{m_i(\lambda)}$ on the
left hand side is one over the number of permutations on $\sum_i im_i$
symbols which commute with a permutation with $m_i$ cycles of length
$i$. This is $\frac{1}{\prod_i i^{m_i} m_i!}$ which agrees with the
coefficient of $\prod_i x_i^{m_i(\lambda)}$ on the right hand
side. \end{proof}

	Combined with Theorem \ref{likeGessel}, this allows us to
prove

\begin{theorem} \label{cauchyidentity} \[ \sum_{\lambda}
S_{\lambda}(\alpha,\beta,\gamma) s_{\lambda}(x) = \sum_{\lambda}
\frac{1}{z_{\lambda}} \tilde{p}_{\lambda} (\alpha,\beta,\gamma)
p_{\lambda}(x) \] \end{theorem}

\begin{proof} Let $n \rightarrow \infty$ in both sides of the
statement of Theorem \ref{likeGessel}. Writing

\[ D_n \left(e^{\gamma z} \prod_{r=1}^{\infty} \frac{1+\beta_r
z}{(1-x_r/z)(1-\alpha_r z)} \right) = D_n \left(e^{\gamma z+\sum_{r
\geq 1} log(1+\beta_r z ) - log(1-\alpha_r z) - log(1-x_r/z)} \right) \]
and applying the strong Szego theorem gives that the determinant
converges to

\[ e^{\sum_{k=1}^{\infty} \frac{1}{k} \tilde{p}_k(\alpha,\beta,\gamma)
p_k(x)}.\] However we know that

\[ \sum_{\lambda} \frac{1}{z_{\lambda}} \tilde{p}_{\lambda}
(\alpha,\beta,\gamma) p_{\lambda}(x) = 
 e^{\sum_{k=1}^{\infty} \frac{1}{k} \tilde{p}_k(\alpha,\beta,\gamma)
p_k(x)}\] because this equation is simply the Polya cycle index (Lemma
\ref{Polya}) with $x_i$ replaced by $\tilde{p_i}(\alpha,\beta,\gamma)
p_i(x)$. \end{proof}

	Richard Stanley has pointed out that Theorems \ref{likeGessel}
and \ref{cauchyidentity} can be deduced from the corresponding results
for usual Schur functions as follows. Since the $h_k$'s are
algebraically independent, there is a unique homomorphism (sending 1
to 1) and $h_k$ to $\tilde{h}_k(\alpha,\beta,\gamma)$ and the
Jacobi-Trudi identity shows that $s_{\lambda}$ maps to
$\tilde{s}_{\lambda}$.

\section{Connections with Work of Baik and Rains} \label{Connect}

	This section connects card shuffling with work of Baik and
Rains \cite{BaRa}. They study ``extended growth models'' indexed by
parameter sets which we call $(\vec{\alpha}^+,\vec{\beta}^+,\gamma^+)$
and $(\vec{\alpha}^-,\vec{\beta}^-,\gamma^-)$. The case relevant to
this paper is $\vec{\alpha}^+=\vec{\beta}^+=\vec{0}$. We assume
without loss of generality (one can simply rescale $\gamma^+$)that
$\gamma^- + \sum \alpha_i^- + \sum \beta_i^-=1$. In this case, which
we call BR$(\gamma^+, \vec{\alpha}^-,\vec{\beta}^-,\gamma^-$), their
model becomes the following:

\begin{enumerate}

\item On $[0,1] \times [0,1]$ choose $Poisson(\gamma^+ \gamma^-)$
i.i.d. uniform points.

\item On $[0,1] \times i$ ($i \in \{1,2,\cdots\}$) choose
$Poisson(\gamma^+ \alpha_i^-)$ i.i.d. uniform points.

\item On $[0,1] \times i$ ($i \in \{-1,-2,\cdots\}$) choose
$Poisson(\gamma^+ \beta_i^-)$ i.i.d. uniform points.

\end{enumerate}

	They define a sequence of points $(x_i,y_i)$ to be increasing
if $x_i \leq x_{i+1},y_i \leq y_{i+1}$ and \[ y_i=y_{i+1}
\Longrightarrow y_i \geq 0.\] They associate to their point process a
random partition $\lambda$ with $\lambda_i$ defined by the property
that \[ \sum_{i=1}^l \lambda_i \] is the size of the longest
subsequence of points which is a union of $l$ increasing subsequences.

	They find a Toeplitz determinant expression for the
probability that $\lambda_1<k$. Theorem \ref{extend} (which is well
known for the case of random permutations (i.e. $\vec{\alpha^-} =
\vec{\beta^-} =0$) relates their point process to card shuffling
measures on permutations and gives a formula for the chance that their
random partition is $\lambda$.

\begin{theorem} \label{extend} 
\begin{enumerate}

\item Consider the random partition arising from the BR($\gamma^+,
\vec{\alpha}^-,\vec{\beta}^-,\gamma^-$) point process. The probability
that this partition is equal to $\lambda$ is the same as the
probability that the RSK shape of a permutation after a
$(\vec{\alpha}^-,\vec{\beta}^-,\gamma^-)$ shuffle on Poisson($\gamma^+)$
symbols is equal to $\lambda$.

\item More explicity, letting $f_{\lambda}$ be the number of standard
Young tableaux of shape $\lambda$, this probability is

\[ \frac{(\gamma^+)^{|\lambda|} f_{\lambda}
S_{\lambda}(\vec{\alpha}^-,\vec{\beta}^-,\gamma^-)}{e^{\gamma^+}
|\lambda|!} .\]

\end{enumerate}
\end{theorem} 	

\begin{proof} We associate to a realization of the BR($\gamma^+,
\vec{\alpha}^-,\vec{\beta}^-,\gamma^-$) point process a random
permutation $\pi$ as follows. First take the deck size to be the
number of points (which has distribution Poisson($\gamma^+)$). Rank
the $y$ coordinates of the points in increasing order, where one
breaks ties for negative $y$ coordiantes by defining the point with
the larger $x$ coordinate to be smaller and breaks ties for positive
$y$ coordinates by defining the point with larger $x$ coordinate to be
larger. Then $\pi(i)$ is defined as the rank of the $y$ coordinate of
the point with the $i$th smallest $x$ coordinate (with probability one
there is no repetition among $x$ coordinates). For example, if the BR
point process yields the points \[
(.2,.3),(.3,.5),(.35,-8),(.4,9),(.45,9),(.5,7),(.6,-2),(.7,-8) \]
then the resulting permutation would be (in 2-line form) \[ 4 \ 5 \ 2
\ 7 \ 8 \ 6 \ 3 \ 1.\] It is easy to see that this distribution on
permutations is the same as that arising from a $(
\vec{\alpha}^-,\vec{\beta}^-,\gamma^-)$ shuffle.
	The second part follows from the first part and Theorem
\ref{probinter}. \end{proof}

	As a corollary, we obtain another proof of a result of Baik
and Rains.

\begin{cor} (\cite{BaRa}) Let $\lambda$ be the partition associated to
the BR($\gamma^+, \vec{\alpha}^-,\vec{\beta}^-,\gamma^-$) point
process. Then the probability that $\lambda$ has largest part at most
$n$ is equal to the Toeplitz determinant \[ \frac{1}{e^{\gamma^+}} D_n
\left(e^{\gamma^+/z} e^{\gamma^-z} \prod_{r=1}^{\infty}
\frac{1+\beta_r^- z}{1-\alpha_r^- z} \right).\] \end{cor}

\begin{proof} By Theorem \ref{extend}, the sought probability is \[
\frac{1}{e^{\gamma^+}} \sum_{\lambda \atop l(\lambda) \leq n}
\frac{(\gamma^+)^{|\lambda|} f_{\lambda}}{|\lambda|!}
S_{\lambda}(\vec{\alpha}^-,\vec{\beta}^-,\gamma^-).\] Writing

\[ \frac{(\gamma^+)^{|\lambda|} f_{\lambda}}{|\lambda|!} = det
\left(\frac{\gamma^+}{(\lambda_i-i+j)!}\right) \] (page 117 of
\cite{Mac}) and

\[ S_{\lambda}(\vec{\alpha}^-,\vec{\beta}^-,\gamma^-) = det
\left(\tilde{h}_{\lambda_i-i+j}(\vec{\alpha}^-,\vec{\beta}^-,\gamma^-)
\right), \] the result now follows by an argument as in Theorem
\ref{likeGessel}. \end{proof}

	As a final result, we note that certain specializations of
Schur measure on partitions (which are of interest to the random
matrix community) have a probability interpretation in terms of
distributions on permutations. For its statement, recall that a
permutation is said to have a descent at position $i$ if
$\pi(i)>\pi(i+1)$; maj denotes the major index of a permutation (sum
of the positions of the descents) and d denotes the number of descents
of a permutation. Finally $\qb{n}{m}_q$ denotes the $q$ binomial
coefficient $\frac{(q^n-1) \cdots (q-1)}{(q^m-1) \cdots (q-1)
(q^{n-m}-1) \cdots (q-1)}$.

\begin{theorem} Consider the probability measure on partitions of size
$n$ which picks $\lambda$ a partition of $n$ with probability \[
\frac{1}{Z_n} s_{\lambda}(1,\frac{1}{p},\cdots,\frac{1}{p^{k-1}})
s_{\lambda}(1,\frac{1}{q}.\cdots,\frac{1}{q^{l-1}})\] (where $Z_n$ is
the normalization constant which can be computed from Cauchy's
identity). This is the pushforward under the RSK correspondence of the
measure on $S_n$ which picks a permutation $\pi$ with probability \[
\frac{1}{Z_n} p^{maj(\pi^{-1})} q^{maj(\pi)} \qb{k-d(\pi^{-1})+n-1}{n}_p
\qb{l-d(\pi)+n-1}{n}_q.\] \end{theorem}

\begin{proof} This follows easily from Proposition 7.9.12 of
\cite{Sta2} together with the fact that if $\pi$ goes to the pair
$(P,Q)$ under the RSK correspondence, the the descent set of $\pi$ is
equal to the descent set of $Q$ and the descent set of $\pi^{-1}$ is
equal to the descent set of $P$ (Lemma 7.23.1 of
\cite{Sta2}). \end{proof}

	Note that for the case of greatest interest ($p=q=1$), the
normalization constant is ${kl+n-1 \choose n}$.

\section{Applications to Card Shuffling: Convergence Rates and Cycle
Index} \label{Shuffling}

	First we derive an upper bound on the convergence rate of
$(\vec{\alpha},\vec{\beta},\gamma)$ shuffles to randomness using
strong uniform times as in \cite{DFP} and then \cite{F0}. The
separation distance between a probability $P(\pi)$ and the uniform
distribution $U(\pi)$ is defined as
$max_{\pi}(1-\frac{Q(\pi)}{U(\pi)})$ and gives an upper bound on total
variation distance. Examples of the upper bound of Theorem
\ref{mybound} are considered later.

\begin{theorem} \label{mybound} The separation distance between $k$
applications of a $(\vec{\alpha},\vec{\beta},\gamma)$ shuffle and
uniform is at most \[ {n \choose 2} \left[ \sum_i (\alpha_i)^2+ \sum_i
(\beta_i)^2 \right]^k .\] Thus $k=2log_{\frac{1}{ \sum_i (\alpha_i)^2+
\sum_i (\beta_i)^2 }}n$ steps suffice to get close to the uniform
distribution.  \end{theorem}

\begin{proof} For each $k$, let $A^k$ be a random $n \times k$ matrix
formed by letting each entry equal $i>0$ with probability $\alpha_i$,
$i<0$ with probability $\beta_i$, and $i=0$ with probability
$\gamma$. Let $T$ be the first time that all rows of $A^k$ containing
no zeros are distinct; from the inverse description of
$(\vec{\alpha},\vec{\beta},\gamma)$ shuffles this is a strong uniform
time in the sense of Sections 4B-4D of Diaconis $\cite{Diac}$, since
if all cards are cut in piles of size one the permutation resulting
after riffling them together is random. The separation distance after
$k$ applications of a $(\vec{\alpha},\vec{\beta},\gamma)$ shuffle is
upper bounded by the probability that $T>k$ \cite{AD}. Let $V_{ij}$ be
the event that rows $i$ and $j$ of $A^k$ are the same and contain no
zeros. The probability that $V_{ij}$ occurs is $\left[ \sum_i
(\alpha_i)^2+ \sum_i (\beta_i)^2 \right]^k$. The result follows
because

\begin{eqnarray*}
Prob(T>k) & = & Prob (\cup_{1 \leq i < j \leq n}) V_{ij}\\
& \leq & \sum_{1 \leq i < j \leq n} Prob(V_{ij})\\
& = & {n \choose 2} \left[ \sum_i (\alpha_i)^2+ \sum_i
(\beta_i)^2 \right]^k
\end{eqnarray*}
\end{proof}

	The results of Section \ref{Extended} and Section \ref{Cauchy}
are used to find a cycle index after
$(\vec{\alpha},\vec{\beta},\gamma)$ shuffles.

	Recall that \[ \tilde{p}_1(\vec{\alpha},\vec{\beta},\gamma) =
\sum_i \alpha_i + \sum_i \beta_i + \gamma =1 \] and (for $n \geq 2$)
\[ \tilde{p}_n(\vec{\alpha},\vec{\beta},\gamma) = \sum_i (\alpha_i)^n
+ (-1)^{n+1} \sum_i (\beta_i)^n. \]

\begin{theorem} \label{cindex}

\begin{enumerate}

\item Let $E_{n,(\vec{\alpha},\vec{\beta},\gamma)}$ denote expected
value after a $(\vec{\alpha},\vec{\beta},\gamma)$ shuffle of an $n$
card deck. Let $N_i(\pi)$ be the number of $i$-cycles of a permutation
$\pi$. Then

\[ \sum_{n \geq 0} u^n E_{n,(\vec{\alpha},\vec{\beta},\gamma)}
(\prod_i x_i^{N_i}) = \prod_{i,j} e^{\frac{(u^ix_i)^j}{ij} \sum_{d|i}
\mu(d) \tilde{p}_{jd}(\vec{\alpha},\vec{\beta},\gamma)^{i/d}}.\]

\item Let $E'_{n,(\vec{\alpha},\vec{\beta},\gamma)}$ denote expected
value after a $(\vec{\alpha},\vec{\beta},\gamma)$ shuffle of an $n$
card deck followed by reversing the order of the cards. Then

\[ \sum_{n \geq 0} u^n E'_{n,(\vec{\alpha},\vec{\beta},\gamma)}
(\prod_i x_i^{N_i}) = \sum_{n \geq 0} u^n
E_{n,(\vec{\beta},\vec{\alpha},\gamma)} (\prod_i x_i^{N_i}).\]

\end{enumerate}
\end{theorem}

\begin{proof} Given the results of Section \ref{Extended} and Section
\ref{Cauchy}, the proof of the first part runs along exactly the same
lines as in the proof of Theorem 4 in \cite{FDeal}. The second
assertion follows from the observation that a
$(\vec{\alpha},\vec{\beta},\gamma)$ shuffle followed by reversing the
order of the cards is conjugate (by the longest length element in the
symmetric group) to a $(\vec{\beta},\vec{\alpha},\gamma)$
shuffle. Alternatively, arguing as in the proof of Theorem 5 in
\cite{FDeal}, one sees that the effect of reversing the cards on the
cycle index of a $(\vec{\alpha},\vec{\beta},\gamma)$ shuffle is to
get

\[ \prod_{i,j} e^{\frac{((-u)^ix_i)^j}{ij} \sum_{d|i} \mu(d)
(-\tilde{p}_{jd}(\vec{\alpha},\vec{\beta},\gamma))^{i/d}}.\] \end{proof}

{\bf Example 1} As a first application of Theorem \ref{cindex}, we
derive an expression for the expected number of fixed points,
generalizing the expressions in \cite{DMP},\cite{FDeal}. To get the
generating function for fixed points, one sets $x_2=x_3=\cdots=1$ in
the cycle index. Using the same trick as in \cite{DMP} and
\cite{FDeal}, the generating function simplifies to \[ \frac{1}{1-u}
\frac{e^{ux \gamma}}{e^{ux}} \prod_{i \geq 1} \frac{1-u \alpha_i}{1-ux
\alpha_i} \frac{1+ux \beta_i}{1+u \beta_i}.\] Taking the derivative
with respect to $x$ and the coefficient of $u^n$, one sees that the
expected number of fixed points is \[ \gamma + \sum_{j=1}^n [\sum_i
(\alpha_i)^j + (-1)^{j+1} (\beta_i)^j] .\] This is exactly the sum of
the first $n$ extended power sum functions at the parameters
$(\vec{\alpha},\vec{\beta},\gamma)$.

{\bf Example 2} We suppose that $\vec{\beta}=\vec{0}$ and that
$\alpha_1=\cdots=\alpha_q=\frac{1-\gamma}{q}$. Then the cycle index
simplifies to

\[ \prod_{i \geq 1} \left(\frac{1}{1-x_i
(\frac{u(1-\gamma)}{q})^i}\right)^{\frac{1}{i} \sum_{d|i} \mu(d)
q^{i/d}} \prod_{i \geq 1} e^{\frac{u^ix_i (1-(1-\gamma)^i)}{i}}.\] Of
particular interest is the further specialization $q=1$. Then the
cycle index becomes

\[ \frac{1}{1-x_1u(1-\gamma)} \prod_{i \geq 1} e^{\frac{u^ix_i
(1-(1-\gamma)^i)}{i}}.\]

	Recall that a $(1/2,0,1/2)$ shuffle takes a binomial(n,1/2)
number of cards, thoroughly mixes them, and then riffles them with the
remaining cards. Example 3 on page 140 of \cite{DFP} proves (in
slightly different notation) that the convolution of $k$ (1/2,0,1/2)
shuffles is the same as a $((1/2)^k,0,1-(1/2)^k)$ shuffle. They
conclude (in agreement with Theorem \ref{mybound}) that a
$(1/2,0,1/2)$ shuffle takes $log_2(n)$ steps to be mixed, as compared
to $\frac{3}{2} log_2(n)$ for ordinary riffle shuffles. They also
establish a cut-off phenomenon. From the computation of Example 1 one
sees that the expected number of fixed points also drops and that the
mean mixes twice as fast.

	As another example, consider a $(1-1/n,0,1/n)$
shuffle. Heuristically this is like top to random and \cite{DFP}
proves that the convergence rate is the same ($nlog(n)$ steps), which
agrees with Theorem \ref{mybound}. From page 139 of \cite{DFP},
performing a $(1-1/n,0,1/n)$ shuffle $k$ times is the same as
perfoming a single $((1-1/n)^k,0,1-(1-1/n)^k)$ shuffle. Example 1
gives a formula for the expected number of fixed points; it would be
interesting to derive a cycle index for convolutions of top to random.

	The approach of either \cite{DMP} (method of moments) or
\cite{FDeal} (generating functions) can be used to prove the following
limit result. For its statement, $\mu$ denotes the Moebius function of
elementary number theory. Note that considerable simplifications take
place when $q=1$ (the interesting case) because $\sum_{d|i} \mu(d)$ is
$1$ if $i=1$ and is $0$ otherwise.

\begin{cor} \label{limitbehavior}
\begin{enumerate}

\item Fix $u$ such that $0<u<1$. Choose a random deck size with
probability of getting $n$ equal to $(1-u)u^n$. Let $N_i(\pi)$ be the
number of $i$-cycles of $\pi$ distributed as a
$(\vec{\alpha},\vec{\beta},\gamma)$ shuffle where $\vec{\beta}
=\vec{0}$ and $\alpha_1=\cdots=\alpha_q
=\frac{1-\gamma}{q}$. Then the random variables $N_i$ are independent,
where $N_i$ is the convolution of a Poisson$((u^i(1-(1-\gamma)^i))/i)$
with $\frac{1}{i} \sum_{d|i} \mu(d) q^{i/d}$ many geometrics with
parameter $(\frac{u(1-\gamma)}{q})^i$.

\item Let $N_i(\pi)$ be the number of $i$-cycles of $\pi$ distributed
as a $(\vec{\alpha},\vec{\beta},\gamma)$ shuffle where $\vec{\beta}
=\vec{0}$ and $\alpha_1=\cdots=\alpha_q =\frac{1-\gamma}{q}$. Then as
$n \rightarrow \infty$ the random variables $N_i$ are independent,
where $N_i$ is the convolution of a Poisson$((1-(1-\gamma)^i)/i)$ with
$\frac{1}{i} \sum_{d|i} \mu(d) q^{i/d}$ many geometrics with parameter
$(\frac{1-\gamma}{q})^i$.

\end{enumerate}
\end{cor}

{\bf Example 3} As a final example, consider the case when
$\alpha_1=\cdots=\alpha_q=\beta_1=\cdots=b_q=\frac{1}{2q}$ and all
other parameters are $0$. Theorems \ref{probinter} and \ref{cindex}
imply that the distribution on RSK shape and cycle index is the same
as for the shuffles in Section 5 of \cite{FDeal}, though we do not see
a simple reason why this should be so.

\section{Acknowledgements} This research was supported by an NSF
Postdoctoral Fellowship. We thank Persi Diaconis for comments on this
work.

\end{document}